\documentclass[12pt]{article}
\usepackage{amssymb}
\usepackage{amsmath,amsthm}
\usepackage[latin1]{inputenc}
\usepackage{hyperref}
\usepackage{color}
\usepackage{graphicx}
\hypersetup{colorlinks=true, linkcolor=blue, citecolor=blue,
urlcolor=blue}

 \newtheorem{remark}{Remark}

 \newtheorem{lemma}[remark]{Lemma}
 \newtheorem{theorem}[remark]{Theorem}
 
 \newtheorem{corollary}[remark]{Corollary}

\title{Alliance free and alliance cover sets}

\author{Juan A.
Rodr\'{\i}guez-Vel\'{a}zquez$^{1,}$\footnote{e-mail:\mbox{\tt
juanalberto.rodriguez\@@urv.cat}. Partially supported by the Spanish
Ministry of Education through projects TSI2007-65406-C03-01
``E-AEGIS" and CONSOLIDER CSD2007-00004 ``ARES" and by the Rovira i
Virgili
 University through project 2006AIRE-09
},   Jos\'{e} M. Sigarreta$^{2}$\footnote{e-mail:\mbox{\tt
    josemaria.sigarreta\@@uc3m.es}}, \\ Ismael G.
Yero$^{1,}$\footnote{e-mail: ismael.gonzalez\@@urv.cat} and Sergio
Bermudo$^{3,}$\footnote{e-mail:\mbox{\tt
    sbernav\@@upo.es.} Partially supported by Ministerio de Ciencia y
Tecnolog\'ia, ref. BFM2003-00034 and Junta de Andaluc\'ia, ref.
FQM-260 and ref. P06-FQM-02225.}\\
    \\
$^1${\em Department of Computer Engineering and Mathematics
}\\
Universitat Rovira i Virgili, \\ Av. Pa\"{\i}sos Catalans 26, 43007
Tarragona, Spain.
\\
$^2$ {\it Faculty of Mathematics} \\
Autonomous University of Guerrero
\\
Carlos E. Adame 5, Col. La Garita, Acapulco, Guerrero, M\'{e}xico
\\
$^3${\it Department of Economy, Quantitative Methods and Economic
History} \\
Pablo de Olavide University \\ Carretera de Utrera Km. 1,
41013-Sevilla, Spain }

\begin{document}

\maketitle

\begin{abstract}

A  \emph{defensive} (\emph{offensive}) $k$-\emph{allian\-ce} in
$\Gamma=(V,E)$ is a set $S\subseteq V$ such that   every $v$ in $S$
(in the boundary of $S$) has at least $k$ more neighbors in $S$ than
 it has in $V\setminus S$.  A set $X\subseteq
V$ is \emph{defensive} (\emph{offensive}) $k$-\emph{alliance free,}
if for all defensive (offensive) $k$-alliance $S$, $S\setminus
X\neq\emptyset$, i.e., $X$ does not contain any defensive
(offensive) $k$-alliance as a subset. A set $Y \subseteq V$ is a
\emph{defensive} (\emph{offensive}) $k$-\emph{alliance cover}, if
for all defensive (offensive) $k$-alliance $S$, $S\cap
Y\neq\emptyset$, i.e., $Y$ contains at least one vertex from each
defensive (offensive) $k$-alliance of $\Gamma$.  In this paper we
show several mathematical properties of defensive (offensive)
$k$-alliance free sets and defensive (offensive) $k$-alliance cover
sets, including tight bounds on the  cardinality of  defensive
(offensive) $k$-alliance free (cover) sets.

\end{abstract}



\section{Introduction}

In \cite{alliancesOne}, P. Kristiansen, S. M. Hedetniemi and S. T.
Hedetniemi introduced several types of alliances in graphs,
including defensive and offensive alliances. We are interested in a
generalization of alliances, namely $k$-alliances, given by Shafique
and Dutton  \cite{kdaf}. In this paper we show several mathematical
properties of $k$-alliance free sets and $k$-alliance cover sets.

We begin by stating some notation and terminology. In this paper
$\Gamma=(V,E)$ denotes a simple graph of  order $n$, size $m$,
minimum degree $\delta$ and maximum degree $\Delta$. For a non-empty
subset $S\subseteq V$, and any vertex $v\in V$, we denote by
$N_S(v)$ the set of neighbors $v$ has in $S$: $N_S(v):=\{u\in S:
u\sim v\}$ and $\delta_S(v) = |N_S(v)|$ denotes the degree of $v$ in
$S$.
 %Similarly, we denote by
%$N_{V\setminus S}(v)$ the set of neighbors $v$ has in $V\setminus
%S$: $N_{V\setminus S}(v):=\{u\in V\setminus S: u\sim v\}.$
The complement of the set $S$ in $V$ is denoted by $\overline{S}$.
The boundary of a set $S\subseteq V$ is defined as $\partial
S:=\displaystyle\cup_{v\in S}N_{\overline{S}}(v).$ A nonempty set of
vertices $S\subseteq V$ is called a {\em defensive}
(\emph{offensive}) $k$-{\em allian\-ce} in $\Gamma$ if for every
$v\in S$ ($v\in \partial S$), $\delta_S(v) \ge
\delta_{\overline{S}}(v)+k$. Hereafter, if there is no  restriction
on the values of $k$, we assume  that $k\in \{-\Delta,...,\Delta\}$.
Notice that any vertex subset is an offensive $k$-alliance for $k\in
\{-\Delta, 1-\Delta,2-\Delta\}$.
% As $\partial V=\emptyset$, in this
%paper we will consider that any offensive $k$-alliance is a proper
%subset of $V$.

A set $X\subseteq V$ is {\em defensive} (\emph{offensive}) {\em
$k$-alliance free,}  $k$-daf ($k$-oaf), if for all defensive
(offensive) $k$-alliance $S$, $S\setminus X\neq\emptyset$, i.e., $X$
does not contain any  defensive (offensive) $k$-alliance as a subset
\cite{kdaf,shafique}. A defensive (offensive) $k$-alliance free set
$X$ is \emph{maximal} if for every defensive (offensive)
$k$-alliance free set $Y$, $X\not\subset Y$.
%it is not a proper subset of a defensive (offensive) $k$-alliance free set.
%if for all $v\notin X$, there exists $S\subseteq X$ such that $S\cup
%\{v\}$ is a defensive (offensive) $k$-alliance.
A \emph{maximum} $k$-daf ($k$-oaf) set is a  maximal ($k$-oaf)
$k$-daf set
 of largest cardinality.

A set $Y \subseteq V$ is a \emph{defensive} (\emph{offensive}) {\em
$k$-alliance cover}, $k$-dac ($k$-oac), if for all defensive
(offensive) $k$-alliances $S$, $S\cap Y\neq\emptyset$, i.e., $Y$
contains at least one vertex from each defensive (offensive)
$k$-alliance of $\Gamma$. A $k$-dac ($k$-oac) set $Y$ is
\emph{minimal} if no proper subset of $Y$ is a defensive (offensive)
$k$-alliance cover set. A \emph{minimum} $k$-dac ($k$-oac) set is a
minimal cover set of smallest cardinality.

\begin{remark}
\mbox{}
\begin{itemize}
\item[$(i)$]If $X$ is a minimal $k$-dac $(k$-oac$)$ set then, for all $v\in X$,
there exists a defensive $($offensive$)$ $k$-alliance $S_v$ for
which $S_v \cap X=\{ v \}$.

\item[$(ii)$] If X is a maximal k-daf $($k-oaf$)$ set, then, for all $v\in
\overline{X}$, there exists $S_v\subseteq X$ such that $S_v \cup \
\{v\}$ is a defensive $($offensive$)$ $k$-alliance.
\end{itemize}
\end{remark}

A defensive (offensive) $k$-alliance  is {\em global} if it  is a
dominating set. For short, in the case of a global offensive
$k$-alliance cover (free) set we will write $k$-goac ($k$-goaf).

%%%%%%%% Notación  %%%%%%%%%%
%%%%%%%%%%%%%%%%%%%%%%%%%%%%%
%%%%%%%%%%%%%%%%%%%%%%%%%%%%%

Associated with the characteristic sets defined above we have the
following invariants:

\begin{itemize}
  \item[] $a_k(\Gamma)$:  minimum cardinality of a defensive
$k$-alliance in $\Gamma$.

\item[] $\gamma_{k}(\Gamma)$:  minimum cardinality of a  global
defensive $k$-alliance  in $\Gamma$.

\item[] $\gamma_k^o(\Gamma)$: minimum cardinality of a
global offensive $k$-alliance in $\Gamma$.

\item[] $\phi_k(\Gamma)$:  cardinality of a maximum $k$-daf set in $\Gamma$.
%If a graph $\Gamma$ does not have a defensive $k$-alliance (for some
%$k$), then we assume that $\phi_k(\Gamma)=|V|=n$.

\item[]
$\phi_k^o(\Gamma)$: cardinality of a maximum $k$-oaf  set in
$\Gamma$.

\item[]  $\phi_k^{go}(\Gamma)$:  cardinality of a maximum
$k$-goaf  set in  $\Gamma$.

\item[] $\zeta_k(\Gamma)$: cardinality of a minimum $k$-dac  set in
$\Gamma$. %If a graph $\Gamma$ does not have a defensive $k$-alliance
%(for some $k$), then we assume that $\zeta_k(\Gamma)=0$.

\item[] $\zeta_k^o(\Gamma)$:  cardinality of a minimum $k$-oac
 set in $\Gamma$.

\item[]  $\zeta_k^{go}(\Gamma)$:   cardinality of a minimum
$k$-goac  set in $\Gamma$.

%\item[] $\gamma(\Gamma)$:  minimum cardinality of a dominating set in $\Gamma$.
\end{itemize}

%%%%%%%%%%%%%%%%%%%%%%%%%%%%%%%%%%%%%
%%%%%%%%%%%%%%%%%%%%%%%%%%%%%%%%%%

The following duality between alliance cover and alliance free sets
was shown in \cite{kdaf,shafique}.

\begin{remark} \label{thDual}
 $X$ is a defensive $($offensive$)$ $k$-alliance cover set
if and only if $\overline{X}$ is defensive $($offensive$)$
$k$-alliance free.
\end{remark}

\begin{corollary}
$\phi_k(\Gamma)+\zeta_k(\Gamma)=\phi_k^o(\Gamma)+\zeta_k^o(\Gamma)=n.$
\end{corollary}

\section{Alliance cover and alliance free sets}

We begin by studying the structure of a set  according to the
structure of its complementary set.

\begin{theorem} \label{lemmaDominating}
If $X$ is a  minimal $k$-dac set, then $\overline{X}$  is a
dominating set.
\end{theorem}

\begin{proof}
By Remark \ref{thDual}, if  $X$ is a  minimal $k$-dac set, then
$\overline{X}$ is a  maximal $k$-daf set. Therefore, for all $v \in
X$, there exists $X_v\subseteq \overline{X}$ such that $X_v\cup
\{v\}$ is a defensive $k$-alliance. So, for every $u\in X_v$,
$\delta_{X_v}(u)+\delta_{\{v\}}(u)= \delta_{X_v\cup \{v\}}(u)\geq
\delta_{\overline{X_v\cup
\{v\}}}(u)+k=\delta_{\overline{X_v}}(u)-\delta_{\{v\}}(u) +k.$ On
the other hand, as $X_v$ is not a defensive $k$-alliance, there
exists $w\in X_v$ such that
$\delta_{X_v}(w)<\delta_{\overline{X_v}}(w)+k$. Hence, by the above
inequalities, $\delta_{\overline{X_v}}(w)+k+
\delta_{\{v\}}(w)>\delta_{\overline{X_v}}(w)- \delta_{\{v\}}(w)+k.$
Thus, $2\delta_{\{v\}}(w)>0$ and, as a consequence, $v$ is adjacent
to $w$.
\end{proof}

Notice that there exist minimal $k$-oac sets such that their
complement sets are not dominating sets. For instance we consider
the graph obtained from the cycle graph $C_8$ by adding the edge
$\{v_1,v_3\}$ and the edge $\{v_5,v_7\}$. In this graph the set
$S=\{v_2,v_3,v_5,v_6,v_7\}$ is a minimal $0$-oac but $\bar{S}$ is
not a dominating set.

%Notice that there exist  maximal  $k$-oaf sets which are not
%dominating sets. For instance we consider the graph obtained from
%the cycle graph $C_8$ by adding the edge $\{ v_1, v_3\}$ and the
%edge $\{ v_5, v_7\}$. In this graph the set $\{v_1, v_4, v_8\}$ is a
%maximal $0$-oaf set but it is not a dominating set.

%\begin{corollary}
%$\zeta_k(\Gamma)\le n-\gamma(\Gamma)$.
%\end{corollary}

%As we will show in Corollary \ref{coro2}, the above bound on
%$\zeta_k(\Gamma)$ can be improved.

\begin{theorem}\label{th1}
If $X$ is a  minimal $k$-dac set, then $\overline{X}$ is a global
offensive $k$-alliance.
\end{theorem}

\begin{proof}
If  $X\subset V$ is a  minimal $k$-dac set, then for every $v\in X$
there exists a defensive $k$-alliance $S_v$ such that $S_v\cap
X=\{v\}$. Hence, $\delta_{S_v}(v)\geq \delta_{\overline{S_v}}(v)+k$
and $\delta_{\overline{X}}(v)\geq \delta_{S_v}(v)\geq
\delta_{\overline{S_v}}(v)+k\geq \delta_{X}(v)+k. $ Therefore, for
every  $v\in X$, we have $\delta_{\overline{X}}(v)\geq
\delta_{X}(v)+k $. On the other hand, by Theorem
\ref{lemmaDominating}, $\overline{X}$ is a dominating set. In
consequence, $\overline{X}$ is a global offensive $k$-alliance in
$\Gamma$.
\end{proof}

%Note that, by definition of offensive $k$-alliances,  in Theorem
%\ref{th1} we only consider the cases $2-\Delta \le k\le \Delta$.
%Hereafter we assume that $k$ only takes the values that give sense
%to the results.

\begin{corollary} \label{coro2} %For every $k\in \{3-\Delta ,..., \Delta\}$,
$\phi_k(\Gamma)\ge \gamma_k^o(\Gamma)$ and $\zeta_k(\Gamma)\le
n-\gamma_k^o(\Gamma)$.
\end{corollary}

 Notice that if one vertex $v\in V$ belongs to any offensive
 $k$-alliance,
then $V\setminus \{v\}$ is a  $k$-oaf set. Hence, $\delta(v)< k$.
So, if $ k\le \delta$ and $X$ is a minimal $k$-oac set, then $|X|\ge
2$.

\begin{theorem} \label{th2}
For every $k\in \{2-\Delta,...,\Delta\}$, if $X$  is a minimal
$k$-goac set such that $|X|\ge 2$, then $\overline{X}$ is an
offensive $(k-2)$-alliance. Moreover, if $k\in \{3,...,\Delta\}$,
then $\overline{X}$ is a global offensive $(k-2)$-alliance.
\end{theorem}

\begin{proof}
If $X\subset V$ is a  minimal $k$-goac set,  then for all  $v\in X$
there exists a global offensive $k$-alliance, $S_v$, such that
$S_v\cap X=\{v\}$. Hence, $1+ \delta_{\overline{X}}(u)\geq
\delta_{S_v}(u)\geq \delta_{\overline{S_{v}}}(u)+k\geq
\delta_{X}(u)+k-1,$ for every $u\in \overline{S_{v}}$. As
$X\setminus \{v\} \subset \overline{S_{v}}$, we have $
\delta_{\overline{X}}(u)\geq \delta_{X}(u)+k-2$ for every $u\in
X\setminus \{v\}$. Therefore, $\overline{X}$ is an  offensive
$(k-2)$-alliance. Moreover, if $k>2$, $\overline{X}$ is a dominating
set. So, in such a case, it is a global offensive $(k-2)$-alliance.
\end{proof}

\begin{corollary} For every  $k\in \{3,...,  \delta\}$,
$\phi_k^{go}(\Gamma)\ge \gamma_{k-2}^o(\Gamma)$ and
$\zeta_k^{go}(\Gamma)\le n-\gamma_{k-2}^o(\Gamma)$.
\end{corollary}

%\textcolor{red}{Sabemos que si X es k-oac, entonces $\overline{X}$ es
%k-oaf. Ahora hay que comprobar que pasa con las globales:  sospecho
%que X k-goac NO implica que $\overline{X}$  sea k-goaf, lo de conjunto
%dominante puede fallar. ¿?}

\begin{theorem} \label{th3}
For every $k\in \{1-\Delta,...,\Delta-1\}$,
 \begin{itemize}
\item[\mbox{\rm (i)}] if $X$ is a global offensive $k$-alliance, then
$\overline{X}$ is $(1-k)$-daf;
\item[\mbox{\rm (ii)}] if  $X$  is a defensive
$k$-alliance, then $\overline{X}$ is $(1-k)$-goaf.
\end{itemize}
\end{theorem}

\begin{proof}
(i) If  $X$ is a global offensive $k$-alliance, then for every $v\in
\overline{X}$ we have $\delta_{X}(v)+1-k> \delta_{\overline{X}}(v)
$. Hence, the set $\overline{X}$ is not a defensive
$(1-k)$-alliance. Moreover, if $Y\subset \overline{X}$, then for
every $y\in Y$ we have $\delta_{\overline{Y}}(y)+1-k\ge
\delta_{X}(y)+1-k> \delta_{\overline{X}}(y)\ge \delta_Y(y)$. Thus,
the set $Y$ is not a defensive $(1-k)$-alliance. Therefore,
$\overline{X}$ is a $(1-k)$-daf set.

(ii) If $ X$ is a defensive  $k$-alliance, then for every $v\in X$
we have $ \delta_{\overline{X}}(v)<\delta_{X}(v)+(1-k) $. So,
$\overline{X}$ is not a global offensive  $(1-k)$-alliance.
Moreover, for every $S \subset \overline{X}$ and $v\in X\subset
\overline{S}$ it is satisfied $\delta_{S}(v)\le
\delta_{\overline{X}}(v)< \delta_{X}(v)+(1-k)\leq
\delta_{\overline{S}}(v)+(1-k)$, in consequence, $S$ is not a global
offensive $(1-k)$-alliance.
\end{proof}

\begin{corollary} For every $k\in \{1-\Delta ,... ,   \Delta-1\}$,
 \begin{itemize}
\item[\mbox{\rm (i)}]
$\zeta_{1-k}(\Gamma)\le \gamma_{k}^o(\Gamma)$ and
$\phi_{1-k}(\Gamma)\ge n-\gamma_{k}^o(\Gamma)$;
\item[\mbox{\rm (ii)}] $\zeta_{1-k}^{go}(\Gamma)\le a_{k}(\Gamma)$.
\end{itemize}
\end{corollary}

%If $X \subset V$ is a maximal $k$-daf of $\Gamma=(V,E)$, then
%$\overline{X}$ is $(1-k)$-daf.

% If $X$ is an independent set, then for all
%$v\in X$, $S_{v}=\overline{X}\cup \{v\}$ is a global offensive
%$k$-alliance such that $S_{v}\cap X=\{v\}$.

Notice that all equalities in the above corollaries are attained for
the complete graph of order $n$ where $\phi_{k}(K_{n})=
n-\zeta_{k}(\Gamma)=\gamma_k^o(K_n)=\left\lceil\frac{n+k-1}{2}\right\rceil$
and $\zeta_{1-k}^{go}(\Gamma)=n-\phi_{1-k}^{go}(\Gamma)=
a_{k}(\Gamma)=\left\lceil\frac{n+k+1}{2}\right\rceil$.

%In the following table we collect  the main obtained results
%relating $X$ and $\overline{X}$.

%\begin{center}
%\begin{tabular}[l]{|c|c|c|}
%\hline
% Result & If $X$ is a & then $\overline{X}$ is a \\
% \hline
% Th. \ref{lemmaDominating} & minimal $k$-dac set, & dominating set. \\
%  \hline
%   Th. \ref{th1} & minimal $k$-dac set, & global offensive  $k$-alliance. \\
 % \hline
%Th. \ref{th2} & minimal $k$-goac set ($|X|\ge 2$ and $k> 2$), & global offensive  $(k-2)$-alliance. \\
 %  \hline
  %Th. \ref{th3} & global offensive  $k$-alliance, & $(1-k)$-daf set. \\
 %  \hline
 % Th. \ref{th4} &   defensive  $k$-alliance, & $(1-k)$-goaf set.\\
 %  \hline
 % \end{tabular}
 % \end{center}

As we show in the following table, by combining some of the above
results we can deduce basic properties on alliance free sets and
alliance cover sets. For the restrictions on $k$, see the premises
of the corresponding results.

\begin{center}
\begin{tabular}[l]{|c|c|}
 \hline
 Rem. \ref{thDual} and Th. \ref{lemmaDominating} & Any maximal $k$-daf set is a  dominating set. \\
  \hline
    Rem. \ref{thDual} and Th. \ref{th1} & Any maximal $k$-daf set is a  global offensive  $k$-alliance. \\
  \hline
   Rem. \ref{thDual} and  Th. \ref{th3}   & Any global offensive  $k$-alliance is a  $(1-k)$-dac set.\\
   \hline
  Th. \ref{th1} and Th. \ref{th3} & Any minimal $k$-dac set is
  $(1-k)$-daf. \\
  \hline
   Th. \ref{th2}  and  Th. \ref{th3} &  Any minimal $k$-goac set  of cardinality at least 2
   is
   $(3-k)$-daf.\\
   \hline
  \end{tabular}
  \end{center}

\subsection{Monotony of  $\phi_{k}^{go}(\Gamma)$ and $\phi_{k}(\Gamma)$}
%%%%%%%%%%%%%%%%%%%%%%%%%
%%%%%%%%%%%%%%%%%%%%%%%%
%%%%%%%%%%%%%%%%%%%%%%%%%
\begin{theorem}\label{teo-k-k+2}
If $X$ is a  $k$-goaf set, $k\in\{1,...,\Delta-2\}$, such that
$|X|\le n-2$, then there exists $v\in \overline{X}$ such that $X\cup
\{v\}$ is a $(k+2)$-goaf set.
\end{theorem}

\begin{proof} %As $k\le \delta$, we have $|\overline{X}|\ge 2$.
 Let us suppose that
for every $x\in \overline{X}$, $X\cup \{x\}$ is not a $(k+2)$-goaf
set. Let  $v\in \overline{X}$ and let $S_v\subset X$, such that
$S_v\cup \{v\}$ is a global offensive $(k+2)$-alliance in $\Gamma$.
Then for every $u\in \overline{S_v\cup
\{v\}}=\overline{S_v}\setminus \{v\}$ we have
$\delta_{S_v}(u)=\delta_{S_v\cup
\{v\}}(u)-\delta_{\{v\}}(u)\ge\delta_{\overline{S_v\cup
\{v\}}}(u)-\delta_{\{v\}}(u)+k+2=
\delta_{\overline{S_v}}(u)-2\delta_{\{v\}}(u)+k+2\ge
\delta_{\overline{S_v}}(u)+k.$ So, for every $u\in
\overline{X}\setminus \{v\}\subset \overline{S_v}\setminus \{v\}$,
$\delta_X(u)\ge \delta_{S_v}(u)\ge \delta_{\overline{S_v}}(u)+k\ge
\delta_{\overline{X}}(u)+k.$ Now we take a vertex $w\in
\overline{X}\setminus \{v\}$ and by the above procedure, taking the
vertex $w$ instead of $v$, we obtain that $\delta_X(v)\ge
\delta_{\overline{X}}(v)+k.$ So, $X$ is a global offensive
$k$-alliance, a contradiction.
%Thus, if $\delta_X(v)\ge \delta_{\overline{X}}(v)+k$,  we obtain
%that $X$ is a global offensive $k$-alliance, a contradiction. If for
%every $v\in \overline{X}$, $\delta_X(v)<
%\delta_{\overline{X}}(v)+k$, then we consider two cases:
 %{\bf Case 1:} Suppose that the subgraph induced by $
%\overline{X}$ is  a complete graph.  In this case, for every $u\in
%\overline{X}\setminus \{v\}$ we have, $\delta_X(u)+1\ge
%\delta_{S_v}(u)+1=\delta_{S_v\cup \{v\}}(u)\ge
%\delta_{\overline{S_v\cup \{v\}}}(u)+k+2=
%\delta_{\overline{S_v}\setminus
%\{v\}}(u)+k+2=\delta_{\overline{S_v}}(u)+k+1\ge
%\delta_{\overline{X}}(u)+k+1.$ Thus,
%$\delta_X(u)\ge\delta_{\overline{X}}(u)+k$, a contradiction.
 %{\bf Case 2:} Suppose that the subgraph induced
%by $ \overline{X}$ is not a  complete graph and let $u,w\in \bar{X}$
%be two non adjacent vertices. As $X$ is a maximal $k$-goaf set,
%there exists $S_u\subset X$, such that $S_u\cup \{u\}$ is a global
%offensive $k$-alliance. Hence, $\delta_X(w)\ge
%\delta_{S_u}(w)=\delta_{S_u\cup \{u\}}(w)\ge
%\delta_{\overline{S_u\cup \{u\}}}(w)+k\ge
%\delta_{\overline{S_u}}(w)+k\ge \delta_{\overline{X}}(w)+k,$ a
%contradiction.
%So, there exists $v\in \overline{X}$, such that $X\cup \{v\}$ is a
%$(k+2)$-goaf set.
\end{proof}

If $X$ is a $k$-goaf for $k\le \delta$, then $|X|\le n- 2$, as a
consequence, the above result can be simplified as follows.
\begin{corollary}
If $X$ is a  $k$-goaf set, $k\in\{1,...,\delta\}$, then there exists
$v\in \overline{X}$ such that $X\cup \{v\}$ is a $(k+2)$-goaf set.
\end{corollary}
It is easy to check  the monotony of $\phi_{k}^{go}$,  i.e.,
$\phi_{k}^{go}(\Gamma)\le \phi_{k+1}^{go}(\Gamma)$. As we can see
below, Theorem \ref{teo-k-k+2} leads to an interesting property
about the monotony of $\phi_{k}^{go}$.
\begin{corollary}
For every $k\in \{1,...,\min\{\delta,\Delta-2\}\}$ and $r\in
\left\{1,...,\lfloor\frac{\Delta-k}{2}\rfloor\right\}$,
$\phi_{k}^{go}(\Gamma)+r\le \phi_{k+2r}^{go}(\Gamma)$.
\end{corollary}

%%%%%%%%%%%%%%%%%%%%%%%
%%%%%%%%%%%%%%%%%%%%%%%%
\begin{theorem}
If $X$ is a $k$-daf set and  $v\in \overline{X}$, then $X\cup \{v\}$
is  $(k+2)-daf$.
\end{theorem}

\begin{proof}
Let us suppose that there exists a  defensive $(k+2)$-alliance $A$
such that $A\subseteq X\cup \{v\}$. If $v\notin A$, then $A\subset
X$, a contradiction because every defensive (k+2)-alliance is a
defensive $k$-alliance. If $v\in A$, let $B=A\setminus \{v\}$. As
for every $u\in B$, $\delta_B(u)=\delta_A(u)-\delta_{\{v\}}(u)$ and
$\delta_{\overline{B}}(u)=\delta_{\overline{A}}(u)+\delta_{\{v\}}(u)$,
we have, $ \delta_A(u)\ge \delta_{\overline{A}}(u)+k+2
\delta_B(u)+\delta_{\{v\}}(u)\ge
\delta_{\overline{B}}(u)-\delta_{\{v\}}(u)+k+2 \delta_B(u)\ge
\delta_{\overline{B}}(u)+k. $
 So, $B\subseteq X$ is a defensive $k$-alliance, a contradiction.
\end{proof}

\begin{corollary}
 For every $k\in \{-\Delta,...,\Delta-2\}$ and $r\in
\left\{1,...,\lfloor\frac{\Delta-k}{2}\rfloor\right\}$,
$\phi_{k}(\Gamma)+r\le \phi_{k+2r}(\Gamma)$.
\end{corollary}

\section{Tight bounds}

A dominating set $S\subset V$ is a \emph{global boundary offensive
$k$-alliance} if for every $v\in \overline{S}$,
$\delta_S(v)=\delta_{\overline{S}}(v)+k$ \cite{Offboundary}.

\begin{lemma}\label{lema-fronteras}
If $\{X, Y\}$ is a vertex partition of a graph $\Gamma$ into two
global boundary offensive $0$-alliances, then $X$ and $Y$ are
minimal global offensive $0$-alliances in $\Gamma$.
\end{lemma}

\begin{proof}
Let us suppose, for instance, that $X$ is not a minimal global
offensive $0$-alliances, then, there exists $A\subset X$, such that,
$X\setminus A\ne \emptyset$ and $A$ is a global offensive
$0$-alliance. Thus, for every $v\in \overline{A}$,
$\delta_X(v)\ge\delta_A(v)\ge\delta_{\overline{A}}(v)\ge
\delta_Y(v)$.

As $Y\subset \overline{A}$ and $\{X, Y\}$ is a vertex partition of
the graph into two global boundary offensive $0$-alliances, then for
every $v\in Y$,
$\delta_Y(v)=\delta_X(v)\ge\delta_A(v)\ge\delta_{\overline{A}}(v)\ge
\delta_Y(v)$.

Therefore, as $Y$ is a dominating set, the above expression carry
out just in the case that $A=X$, a contradiction. So, $X$ and $Y$
are minimal global offensive $0$-alliances.
\end{proof}

\begin{theorem}
For  every  $k\in \{0,...,\Delta\}$, $\phi_k^{go}(\Gamma)\ge
\lfloor\frac{n}{2}\rfloor+\lfloor\frac{k}{2}\rfloor-1$.
\end{theorem}

\begin{proof}
First, we will prove the case $k=0$. Let $\{X, Y\}$ be a partition
of the vertex set, such that $|X|=\lfloor\frac{n}{2}\rfloor$,
$|Y|=\lceil\frac{n}{2}\rceil$ and there is a minimum number of edges
between $X$ and $Y$. If $X$ (or $Y$) is a $0$-goaf set, then
$\phi_0^{go}(\Gamma)\ge \lfloor\frac{n}{2}\rfloor-1$. We suppose
there exist $A\subset X$ and $B\subset Y$, such that $A$ and $B$ are
global offensive $0$-alliances. Hence $\delta_X(v)\ge
\delta_{A}(v)\ge \delta_{\bar{A}}(v)\ge \delta_Y(v), \;\; \forall
v\in \bar{A},$ and $\delta_Y(v)\ge \delta_{B}(v)\ge
\delta_{\bar{B}}(v)\ge \delta_X(v), \;\; \forall v\in \bar{B}.$ As
$Y\subset \bar{A}$ and $X\subset \bar{B}$ we have, for every $v\in
Y$, $\delta_X(v)\ge \delta_Y(v)$ and for every $v\in X$,
$\delta_Y(v)\ge \delta_X(v)$.

For any $y\in Y$ and $x\in X$, let us take $X'=X\setminus \{x\}\cup
\{y\}$ and $Y'=Y\setminus \{y\}\cup \{x\}$. If
$\delta_X(y)>\delta_Y(y)$ or $\delta_Y(x)>\delta_X(x)$ then, the
edge cutset between $X'$ and $Y'$ is lesser than the other one
between $X$ and $Y$, a contradiction. Therefore
$\delta_X(y)=\delta_Y(y)$ and $\delta_Y(x)=\delta_X(x)$ and, as a
consequence, $\{X, Y\}$ is a partition of the vertex set into two
global boundary offensive $0$-alliances. Now, by using Lemma
\ref{lema-fronteras} we obtain that $X$ and $Y$ are minimal global
offensive $0$-alliances. As a consequence, $\phi_0^{go}(\Gamma)\ge
\lfloor\frac{n}{2}\rfloor-1$.

Now, let us prove the case $k>0$. { Case 1}: $\phi_k^{go}(\Gamma)\ge
n-2$. Since $\ n-1\ge \lfloor\frac{2n-1}{2}\rfloor\ge
\lfloor\frac{n+\Delta}{2}\rfloor\ge \lfloor\frac{n+k}{2}\rfloor\ge
\lfloor\frac{n}{2}\rfloor+\lfloor\frac{k}{2}\rfloor,$ we have
$\phi_k^{go}(\Gamma)\ge
\lfloor\frac{n}{2}\rfloor+\lfloor\frac{k}{2}\rfloor-1$. {Case 2:
$\phi_k^{go}(\Gamma)< n-2$.}  As every $k$-goaf set is also a
$(k+1)$-goaf set, $\phi_1^{go}(\Gamma)\ge \phi_0^{go}(\Gamma)\ge
\lfloor\frac{n}{2}\rfloor+\lfloor\frac{1}{2}\rfloor-1$, then the
statement is  true for $k=1$. Hence, we will proceed by induction on
$k$.  Let us assume that the statement is true  for  an arbitrary
$k\in \{2,...,\Delta-2\}$, that is, there exists a   maximal
$k$-goaf set  $X$ in $\Gamma$ such that, $|X|=\phi_k^{go}(\Gamma)\ge
\lfloor\frac{n}{2}\rfloor+\lfloor\frac{k}{2}\rfloor-1$. Now, by
 Theorem \ref{teo-k-k+2}, there exists $v\in
\overline{X}$, such that $X\cup \{v\}$ is a $(k+2)$-goaf set.
Therefore, $\phi_{k+2}^{go}(\Gamma)\ge |X\cup\{v\}|\ge
\lfloor\frac{n}{2}\rfloor+\lfloor\frac{k}{2}\rfloor=
\lfloor\frac{n}{2}\rfloor+\lfloor\frac{k+2}{2}\rfloor-1$. So, the
proof is complete.
\end{proof}

The above bound is attained, for instance, in the case of the
complete graph if $n$ and $k$ are both even or if $n$ and $k$ have
different parity: $ \phi_{k}^{go}(K_{n})=
\left\lfloor\frac{n+k-2}{2}\right\rfloor$.

%From (\ref{no-k-daf}) we deduce that if $X$ is a  $k$-daf set
%whose induced subgraph is complete, then $|X| \leq \left\lfloor
%\displaystyle\frac{\Delta +k+1}{2}\right\rfloor.$

\begin{theorem} \label{Th-oaf}
$\left\lceil\frac{\delta+k-2}{2}\right\rceil \le
\phi_{k}^{o}(\Gamma)\leq
\left\lfloor\frac{2n-\delta+k-3}{2}\right\rfloor.$
\end{theorem}

\begin{proof}
If $X$ is a  $k$-oaf set, then $\delta_X(v)+1\leq
\delta_{\overline{X}}(v)+k,$ for some $v\in \partial X$. Therefore,
$ \delta(v)+1-k= \delta_X(v)+ \delta_{\overline{X}}(v)+1-k\leq
2\delta_{\overline{X}}(v)\leq 2(n-|X|-1). $ Thus, the upper bound is
deduced.

If $X$ is a maximal $k$-oaf set, then $\overline{X}$ is a minimal
$k$-oac set. Thus, for all $v\in \overline{X}$, there exists an
offensive $k$-alliance $S_{v}$ such that $S_{v}\cap
\overline{X}=\{v\}$. Hence, $\delta_{S_{v}}(u)\geq
\delta_{\overline{S_{v}}}(u)+k$, for every $u\in
\partial S_{v}$. Therefore,
%\begin{equation}\label{29}
 $\delta(u)+k\le 2\delta_{S_{v}}(u)\leq 2|S_v|\leq 2(|X|+1).$
%\end{equation}
Thus, the lower bound follows.
\end{proof}

The above bounds are attained, for instance, for the complete graph:
$ \phi_{k}^o(K_{n})= \left\lceil\frac{n+k-3}{2}\right\rceil$.%,  $-n+3< k <n $.

For every $k\in \{0,...,\Delta\}$ it was established in
\cite{shafique} that
  $\phi_k(\Gamma)\ge
\lfloor\frac{n}{2}\rfloor+\lfloor\frac{k}{2}\rfloor$. The next
result shows other bounds on $\phi_k(\Gamma)$.

\begin{theorem} \label{cotaconnectivity}
 For any connected graph  $\Gamma$,
$\left\lceil\frac{n(k+\mu)-\mu}{n+\mu}\right\rceil
\le\phi_{k}(\Gamma)\le  \left\lfloor\frac{2n+k-\delta-1}{2}
\right\rfloor,$ where $\mu$ denotes the algebraic connectivity of
$\Gamma$.
\end{theorem}

\begin{proof}
It was shown in \cite{RoGoSi} that the defensive $k$-alliance number
is bounded by $ a_k(\Gamma)\ge
\left\lceil\frac{n(\mu+k+1)}{n+\mu}\right\rceil$. On the other hand,
if
 $S$ is a  defensive $k$-alliance of cardinality $a_k(\Gamma)$, then
 for all
$v\in S$ we have that  $S\setminus \{v\}$ is a $k$-daf set. Thus,
$\phi_k(\Gamma)\ge a_k(\Gamma)-1$. Hence, the lower bound on
$\phi_k(\Gamma)$  follows.

Moreover, if $X$ is a $k$-daf set, then $\delta_{ X}(v)+1 \le
\delta_{\overline{X}}(v)+k, \quad \mbox{\rm  for some {  } }v\in X.$
Therefore, $\delta(v)+1-k= \delta_{ X}(v)
+\delta_{\overline{X}}(v)+1 -k\le 2 \delta_{\overline{X}}(v)\le
2(n-|X|).$ Thus, the upper bound follows.
\end{proof}

The above bound is sharp as we can check, for instance, for the
complete graph $\Gamma=K_n$.  As the algebraic connectivity of $K_n$
is $\mu=n$, the above theorem gives the exact value of
$\phi_{k}(K_{n})= \left\lceil\frac{n+k-1}{2}\right\rceil.$

%Notice that bounds on $\zeta_k ^o(\Gamma)$ and $\zeta_k(\Gamma)$ can
%be deduced from Theorem \ref{Th-oaf}  and Theorem
%\ref{cotaconnectivity}, respectively.

\begin{theorem}\label{thOf} % FeRoSi
For any connected graph  $\Gamma$, $\zeta_k(\Gamma)\le
\frac{n}{\mu_*} \left( \mu_*-\left\lceil\frac{\delta
+k}{2}\right\rceil\right),$ where $\mu_*$ denotes the Laplacian
spectral radius of $\Gamma$.
\end{theorem}

\begin{proof} The result immediately follows from Corollary \ref{coro2} and the
following bound obtained in \cite{FeRoSi}: $\gamma_k^o(\Gamma)\ge
\displaystyle\frac{n}{\mu_*}\left\lceil\frac{\delta
+k}{2}\right\rceil.$
\end{proof}

The above bound is tight as we can check, for instance, for the
complete graph $\Gamma=K_n$.  As the Laplacian spectral radius of
$K_n$ is $\mu_*=n$, the above theorem gives the exact value of
$\zeta_{k}(K_{n})= \left\lceil\frac{n-k}{2}\right\rceil.$


\begin{thebibliography}{99}

\bibitem{FeRoSi}
H. Fernau, J.A. Rodr\'{\i}guez and J.M. Sigarreta, Offensive
k-alliances in graphs. \emph{ Discrete Applied Mathematics}
\textbf{157} (2) (2009), 177-182.



\bibitem{alliancesOne} S. M. Hedetniemi, S. T. Hedetniemi, and P. Kristiansen,
Alliances in graphs. {\it J. Combin. Math. Combin. Comput.} {\bf 48}
(2004), 157-177.

%\bibitem{GlobalalliancesOne} T. W. Haynes,   S. T. Hedetniemi  and M. A. Henning,
%Global defensive alliances in graphs, {\it Electron. J. Combin.}
%{\bf 10} (2003), Research Paper 47, 13 pp.

\bibitem{RoGoSi}
J.A. Rodr\'{\i}guez, I. G. Yero and J.M. Sigarreta, Defensive
$k$-alliances in graphs. \emph{Applied Mathematics Letter} {\bf 22}
(2009), 96-100.


\bibitem{kdaf} K. H. Shafique and   R. D. Dutton,  Maximum alliance-free and
minimum alliance-cover sets, {\it Congr. Numer. }{\bf 162}
 (2003),139-146.

\bibitem{shafique} K. H. Shafique y R. Dutton, A tight bound on the cardinalities of maximun alliance-free
and minimun alliance-cover sets, {\it J. Combin. Math. Combin.
Comput.} {\bf 56} (2006) 139-145.

%\bibitem{espectral} J. A. Rodr\'{\i}guez  and J. M. Sigarreta,  Spectral study of alliances in graphs.
%\emph{Discussiones Mathematicae Graph Theory} \textbf{27} (1) (2007)
%143-157.





%\bibitem{planar} J. A. Rodr\'{\i}guez-Vel\'{a}zquez and J. M. Sigarreta, Global alliances in planar
%graphs. AKCE-\emph{International Journal of Graphs and
%Combinatorics} \textbf{4} (1) (2007) 83-98.

\bibitem{Offboundary} I. G. Yero and J. A.
Rodr\'{\i}guez-Vel\'{a}zquez, Boundary offensive $k$-alliances in
graphs. Submitted, 2008.


\end{thebibliography}
\end{document}